\documentclass[11pt]{article}
\usepackage{latexsym,amssymb,amsfonts}
\usepackage{mathrsfs}
\usepackage{graphicx}
\usepackage{psfrag}
\usepackage{amsmath, amsbsy}
\usepackage{amsopn, amstext}
\usepackage{amsmath,multicol,amsopn}
\usepackage{latexsym,amssymb}

\usepackage{anysize,hyperref}

\def\sqr#1#2{{\vcenter{\vbox{\hrule height.#2pt
				\hbox{\vrule width.#2pt height#1pt \kern#1pt \vrule width.#2pt}
				\hrule height.#2pt}}}}
\def\signed #1{{\unskip\nobreak\hfil\penalty50
		\hskip2em\hbox{}\nobreak\hfil#1
		\parfillskip=0pt \finalhyphendemerits=0 \par}}
\def\endpf{\signed {$\sqr69$}}

\def\dbR{{\mathop{\rm l\negthinspace R}}}

\def\3n{\negthinspace \negthinspace \negthinspace }
\def\2n{\negthinspace \negthinspace }
\def\1n{\negthinspace }

\def\dbE{{\mathbb{E}}}
\def\dbF{{\mathbb{F}}}

\def\dbN{{\mathbb{N}}}

\def\dbP{{\mathbb{P}}}

\def\dbR{{\mathbb{R}}}

\def\={\buildrel \triangle \over =}

\def\a{\alpha}

\def\d{\delta}
\def\e{\varepsilon}

\def\l{\lambda}

\def\n{\nabla}
\def\si{\sigma}

\def\f{\varphi}
\def\th{\theta}
\def\om{\omega}

\def\ns{\noalign{\ss} }

\def\pa{\partial}
\def\Om{\Omega}

\def\cA{{\cal A}}
\def\cB{{\cal B}}
\def\cC{{\cal C}}

\def\cE{{\cal E}}

\def\cH{{\cal H}}

\def\cJ{{\cal J}}

\def\mE{\mathbb{E}}
\def\no{\noindent}

\def\ss{\smallskip}
\def\ms{\medskip}
\def\bs{\bigskip}
\def\q{\quad}
\def\qq{\qquad}
\def\hb{\hbox}

\def\max{\mathop{\rm max}}
\def\min{\mathop{\rm min}}
\def\exp{\mathop{\rm exp}}

\def\pa{\partial}

\def\cd{\cdot}

\def\as{\hbox{\rm a.s.{ }}}

\def\|{\Big |}
\def\({\Big (}
\def\){\Big )}
\def\[{\Big[}
\def\]{\Big]}
\def\be{\begin{equation}}
	\def\bel{\begin{equation}\label}
		\def\ee{\end{equation}}
	\def\bt{\begin{theorem}}
		\def\bcd{\begin{condition}}
			\def\ecd{\end{condition}}
		\def\et{\end{theorem}}
	\def\bc{\begin{corollary}}
		\def\ec{\end{corollary}}
	\def\bde{\begin{definition}}
		\def\ede{\end{definition}}
	\def\bl{\begin{lemma}}
		\def\el{\end{lemma}}
	\def\bp{\begin{proposition}}
		\def\ep{\end{proposition}}
	\def\br{\begin{remark}}
		\def\er{\end{remark}}
	\def\ba{\begin{array}}
		\def\ea{\end{array}}
	\def\ed{\end{document}}
\def\ns{\noalign{\ms}}
\def\ds{\displaystyle}

\def\square#1{\vbox{\hrule\hbox{\vrule height#1%
			\kern#1\vrule}\hrule}}
\def\rectangle#1#2{\vbox{\hrule\hbox{\vrule height#1%
			\kern#2\vrule}\hrule}}

\font\tenbb=msbm10 \font\sevenbb=msbm7
\font\fivebb=msbm5

\newfam\bbfam
\scriptscriptfont\bbfam=\fivebb
\textfont\bbfam=\tenbb
\scriptfont\bbfam=\sevenbb

\newtheorem{lemma}{Lemma}[section]
\newtheorem{remark}{Remark}[section]

\newtheorem{theorem}{Theorem}[section]
\newtheorem{corollary}{Corollary}[section]

\newtheorem{definition}{Definition}[section]
\newtheorem{proposition}{Proposition}[section]
\newtheorem{condition}{Condition}[section]

\allowdisplaybreaks

\makeatletter

\@addtoreset{equation}{section}
\makeatother

\begin{document}
	\title{\bf Strong Unique Continuation Property for Stochastic Parabolic Equations\thanks{This work is partially supported by the NSF of China
			under grants 12025105 and 11931011, and by the Science Development Project of Sichuan University under
			grant 2020SCUNL201. }}
	
	\author{  Zhonghua Liao\thanks{School of Mathematics,  Sichuan University,
			Chengdu 610064,  China. E-mail address: zhonghualiao@yeah.net} \  and\  Qi L\"u\thanks{School
			of Mathematics, Sichuan University, Chengdu
			610064, China. E-mail address: lu@scu.edu.cn.}}
	
	\date{}

	\maketitle

	\begin{abstract}\no
		We establish a  strong unique continuation property for stochastic
		parabolic equations. Our method is based on a new stochastic version
		of Carleman estimate. As far as we know, this is the first result
		for strong unique continuation property of stochastic partial
		differential equations.
	\end{abstract}
	
	\bs
	
	\no{\bf 2010 Mathematics Subject
		Classification}. 60H15. \bs
	
	\no{\bf Key Words}. strong unique continuation
	property, stochastic parabolic equations,
	Carleman estimate.
	
	\ms
	
	
	\section{Introduction }
	
	\q Let $T > 0$, and $G \in \mathbb{R}^{n}$ ($n \in \mathbb{N}$) be a
	given domain.  Denote $Q = (0, T)\times G$. Assume that $a^{jk}\in
	C^{1,2}([0,T]\times G)$ satisfy $a^{jk}=a^{kj}$ ($j,k=1,2,\cdots,n$)
	, and for any open subset $G_1$ of $G$, there is a constant
	$s_0=s_0(G_1)>0$ such that
	\begin{equation}\label{1.2-eq1}
		\sum_{j,k=1}^n a^{jk}(t,x)\xi_j\xi_k \geq
		s_0|\xi|_{\dbR^n}^2,\qq \forall\
		(t,x,\xi^{1},\cdots, \xi^{n})\in [0,T]\times G_1
		\times\dbR^n.
	\end{equation}

	Let $(\Om, {\cal F}, \dbF, \dbP)$ with $\dbF\=\{{\cal F}_t\}_{t \geq
		0}$ be a complete filtered probability space on which a  one
	dimensional standard Brownian motion $\{W(t)\}_{t\geq 0}$ is
	defined. Assume that $H$ is a Fr\'echet space. Denote by
	$L^2_{\dbF}(0, T; H)$ be the Fr\'echet space consisting all
	$H$-valued $\dbF$-adapted process $X(\cdot)$ such that $\mE
	|X(\cdot)|_{L^2_{\dbF}(0, T; H)}^2 < +\infty$, by
	$L^{\infty}_{\dbF}(0, T; H)$ the Fr\'{e}chet space of all $H$-valued
	$\dbF$-adapted bounded processes  and by $L^2_{\dbF}(\Omega; C([0,
	T]; H))$ the Fr\'{e}chet space of all $H$-valued $\dbF$-adapted
	continuous precesses $X(\cdot)$ with $\mE |X(\cdot)|^2_{C([0, T];
		H)} < +\infty$. All the above spaces are equipped with the canonical
	quasi-norms.

	Consider the following stochastic parabolic equation:
	\begin{equation}\label{pa-sucp0}
		dy-\sum_{j,k=1}^n (a^{jk}y_{x_j})_{x_k} dt=a\cd\nabla ydt + bydt +cydW(t)\qq \hb{ in } Q,\\
	\end{equation}
	where $a\in L^{\infty}_{\dbF}(0,T;L_{loc}^{\infty}(G;\dbR^n))$,
	$b\in L^{\infty}_{\dbF}(0,T;L_{loc}^{\infty}(G))$ and $c\in
	L^{\infty}_{\dbF}(0,T; W_{loc}^{1,\infty}(G))$.
	
	In this paper, for simplicity, we use the
	notation $y_{x_j}\equiv y_{x_j}(x)={{\pa
			y(x)}/{\pa x_j}}$, where $x_j$ is the $j$-th
	coordinate of a generic point $x=(x_1,\cdots,x_n)$ in $\dbR^n$.
	Similarly, we use $ z_{x_j}$, $v_{x_j}$, etc. for the partial
	derivatives of $z$ and $v$ with respect to $x_j$. Also, we use
	$C=((a^{jk})_{1\le j,k\le n}, Q, a, b,c)$ to denote a generic
	positive constant independent of the solution $y$, which may change
	from line to line.

	\vspace{0.1cm}

	To begin with, we recall the definition of the
	solution to \eqref{pa-sucp0}.
	
	\medskip
	
	\begin{definition}
		We call  $y\in
		L_{\dbF}^2(\Om;C([0,T];L_{loc}^2(G)))\cap
		L_{\dbF}^2(0,T;H_{loc}^1(G))$ a solution to
		\eqref{pa-sucp0} if for any $t\in[0,T]$, any
		nonempty open subset $G'$ of $G$ and any
		$\eta\in H_0^1(G')$, it holds
		\begin{equation}\label{zx2}
			\begin{array}{ll}\ds
				\q\int_{G'}y(t,x)\eta(x)dx-\int_{G'}y(0,x)\eta(x)dx\\
				\ns  =\!\ds\int_0^t\int_{G'} \Big
				\{\!-\!\sum_{j,k=1}^n
				a^{jk}(s,x)y_{x_j}(s,x)\eta_{x_k}(x) +[
				a(s,x)\cd
				\n y(s,x) +b(s,x)y(s,x)]\eta(x) \Big \}dx ds\\
				\ns \q\ds+\int_0^t\int_{G'}c(s,x)y(s,x)\eta(x)
				dx dW(s),\qq \dbP\mbox{-}\as
			\end{array}
		\end{equation}
	\end{definition}

	In this paper, we study the strong unique continuation property
	(SCUP for short) for solutions to \eqref{pa-sucp0}. For a given
	deterministic/stochastic PDE, SUCP means, roughly speaking, if a
	solution to the equation vanishes to infinite order at a point of a
	connected open set, then it must vanish identically in that set.
	SUCP is one of the most fundamental aspects  for deterministic PDEs.
	It is studied extensively in the literature. Classical results are
	Cauchy-Kovalevskaya theorem and Holmgren's uniqueness theorem. These
	results need to assume that the coefficients of the PDE to be
	analytic. In 1939, T. Carleman introduced in the seminal paper
	\cite{Carleman1} a new method  to prove SUCP for two dimensional
	elliptic equations with $L^\infty$ coefficients. This landmark work
	indicates that  a non-analytic solution of an elliptic equation can
	behave  in an ``analytic" manner in some sense. The technique he
	used, which is called ``Carleman estimate" now, has became a very
	powerful tool in the study of SUCP for elliptic equations (e.g.
	\cite{Aronszajn,Hormander2,JK,Kenig,Koch1,Sogge1}) and parabolic
	equations (e.g.
	\cite{Escauriaza,Escauriaza2,Escauriaza1,Vessella,Vessella1}).
	
	It is worth noting that SUCP is an important problem not only in the
	uniqueness of the solution to a PDE itself, but also in the study of
	other properties of solutions, such as the nodal sets (e.g.
	\cite{DF1,DF}),  the Anderson localization (e.g. \cite{BK1}), etc.
	Furthermore, it can be applied to solve some application problems,
	such as controllability problems (e.g. \cite{Zu1}), optimal control
	problems (e.g. \cite{LY1}), inverse problems (e.g. \cite{Vessella1})
	and so on.
	
	Compared with the deterministic PDEs, as far as we know, there is no
	result concerning SUCP for stochastic PDEs. In our opinion, it would
	be quite interesting to extend the deterministic SUCP results to the
	stochastic ones. Nevertheless, there are many things which remain to
	be done and some of them seem to be challenging.
	
	Before continuing, we give the definition of  SCUP for the solution
	$y$ to \eqref{pa-sucp0}.
	\begin{definition}\label{pa-sucp-def1}
		A solution $y$ to \eqref{pa-sucp0} is said to satisfy the SUCP if
		$y\equiv 0$ in $Q$, $\dbP$-a.s., provided that $y$ vanishes of
		infinite order at $(0,T)\times\{x_0\}$ for some $x_0\in G$, i.e.,
		for any $N\in\dbN$ and $r>0$, there is a $C_N>0$ such that
		$$\mE\int_{(0,T)\times
			\cB(x_0,r)}|y(t,x)|^2dxdt\leq C_N r^{2N}.$$
	\end{definition}

	The main result of this paper is as follows.
	\begin{theorem}\label{pa-sucp-th1}
		Solutions $y$ to \eqref{pa-sucp0} satisfy SUCP.
	\end{theorem}

	In this paper, similar to the deterministic case, we employ a
	Carleman estimate to establish our SUCP result. In recent years,
	motivated by the study of unique continuation problems (NOT the
	strong unique continuation problems), controllability and
	observability problems, and inverse problems, there are considerable
	progresses  concerning the Carleman estimate for stochastic
	parabolic equations (see \cite{FL,Liuxu1,Luqi3,LY,Tang-Zhang1,Zh}).
	Despite such developments, the SUCP  for stochastic parabolic
	equations remains an open area. None of the above Carleman estimates
	can be used to get the SUCP for our equation \eqref{pa-sucp0} due to
	the choice of weight functions. Indeed, weight functions in these
	papers are designed to get some global energy estimate for
	stochastic parabolic equations with boundary conditions. Moreover,
	due to the extra difficulties caused by the stochastic setting, such
	as the requirement of the adaptedness of solutions with respect to
	the filtration $\dbF$, we cannot simply localize the problem as
	usual because the classical localization technique may change the
	adaptedness of solutions. To overcome such hindrances,  we have
	taken inspiration from the ideas not only in \cite{Tang-Zhang1} but
	also in \cite{BK1,Escauriaza,Vessella} to prove Theorem
	\ref{pa-sucp-th1}.

	There are some other methods to establish the SUCP for parabolic
	equations (e.g. \cite{Chen, Lin,Poon}). However, it seems that these
	method cannot be applied to get the SUCP for stochastic parabolic
	equations. For instance, the key step in \cite{Chen} is to recast
	equations in terms of parabolic self-similar variables. However, it
	seems that this cannot be done for stochastic parabolic equations
	since the related changing of variable with respect to $t$ will
	destroy the adaptedness of solutions, which is a key feature in the
	stochastic setting. The method in \cite{Lin} is to reduce the SUCP
	for  parabolic equations with time-independent coefficients  to the
	SUCP for elliptic equations. This reduction relies on a
	representation formula for solutions of parabolic equations in terms
	of eigenfunctions of the corresponding elliptic operator, and
	therefore cannot be applied to more general equations with
	time-dependent coefficients.  The difficulty for employing methods
	in \cite{Poon} to study the SUCP of \eqref{pa-sucp0}  consists in
	the fact that one cannot simply localize the problem and  do
	changing of variables as usual because they may also change the
	adaptedness of solutions.

	The rest of this paper is organized as follows. In Section
	\ref{Sec2}, as a key preliminary, we prove a weighted identity for a
	stochastic parabolic operator. Section \ref{Sec3}  is devoted to
	establishing a Carleman estimate for stochastic parabolic equations.
	At last, in Section \ref{Sec4}, we prove Theorem \ref{pa-sucp-th1}.
	In this paper, in order to present the key idea in the simplest way,
	we do not pursue the full technical generality.


	\section{A weighted identity for a stochastic
		parabolic operator}\label{Sec2}
	
	
	First, we introduce the following  weighted identity for the
	stochastic parabolic operator
	``$dh-\sum_{j,k=1}^n(a^{jk}h_{x_j})_{x_k}dt$".

	\begin{lemma}\label{c1t1}
		Let $f, \ell\in C^{1,3}(Q)$ and $\Psi\in
		C^{1,2}(Q)$. Let $h$ be an $H^2(G)$-valued It\^o
		process. Set $\th=e^{\ell }$ and $v=\th h$.
		Then, for any $t\in [0,T]$ and a.e. $(x,\om)\in
		G\times\Omega$,
		\begin{eqnarray}\label{c1e2a}
			&&
			2 f\th\[-\sum_{j,k=1}^n (a^{jk}v_{x_j})_{x_k}+\cA v\]\[dh-\sum_{j,k=1}^n(a^{jk}h_{x_j})_{x_k}dt\]\nonumber\\
			&&\q+2 \sum_{j,k=1}^n (fa^{jk}v_{x_j}dv)_{x_k}+ 2
			\sum_{j,k=1}^n\Big\{\sum_{j',k'=1}^n\[2fa^{jk}
			a^{j'k'}\ell_{x_{j'}}v_{x_j}v_{x_{k'}}-fa^{jk}a^{j'k'}\ell_{x_j}v_{x_{j'}}v_{x_{k'}}\nonumber
			\\ &&\q +fa^{jk}(a^{j'k'}\ell_{x_{j'}})_{x_{k'}}vv_{x_j}\] -f\cA a^{jk}\ell_{x_j} v^2\Big\}_{x_k}dt - d\(f\cA v^2+f\sum_{j,k=1}^n
			a^{jk}v_{x_j}v_{x_k}\) \\
			&& =2 \sum_{j,k=1}^n c^{jk}v_{x_j}v_{x_k}dt
			+ \cB v^2dt  +2 f\[-\sum_{j,k=1}^n \big(a^{jk}v_{x_j}\big)_{x_k}+\cA v\]^2dt\nonumber\\
			&&\q - f\th^2\sum_{j,k=1}^n a^{jk}( dh_{x_j} +
			\ell_{x_j}dh)(dh_{x_k}+\ell_{x_k}dh)- f\th^2\cA
			(dh)^2+2\sum_{j,k=1}^n a^{jk}f_{x_j}v_{x_k}dv \nonumber\\
			&&\q +2\sum_{j,k,j'k'=1}^n\big[
			a^{jk}f(a^{j'k'}\ell_{x_{j'}})_{x_{k'}}\big]_{x_j}v_{x_k}vdt,\nonumber
		\end{eqnarray}
		where
		\begin{equation}\label{c1e15}
			\left\{
			\begin{array}{ll}
				\ds\cA =-\sum_{j,k=1}^n
				a^{jk}\ell_{x_j}\ell_{x_k}-\ell_t,\\
				\ns \ds \cB=-2\sum_{j,k=1}^n a^{jk}(f\cA)_{x_j}\ell_{x_k}-(f\cA)_t,\\
				\ns \ds c^{jk}=\sum_{j'k'=1}^n\[ 2a^{jk'}(fa^{j
					'k}\ell_{x_{j'}})_{x_{k'}}-(fa^{jk})_{x_{k'}}a^{j'k'}\ell_{x_{j'}}-\frac{1}{2}(fa^{jk})_t\]
				.
			\end{array}
			\right.
		\end{equation}
	\end{lemma}
	\begin{remark}
		Compared with the widely used weighted identity for stochastic
		parabolic operator (\cite[Theorem 2.1]{Tang-Zhang1}), we introduce
		an auxiliary function $f$ in the pointwise identity, which plays a
		key role in the proof of Theorem \ref{pa-sucp-th1}.
	\end{remark}

	\par {\it Proof  of Lemma \ref{c1t1}}. We divide the proof into
	three steps.
	
	\ss
	
	{\bf Step 1.} Recalling $\th=e^{\ell}$ and $v=\th h$, one has
	$dh=\th^{-1}(dv-\ell_tv)dt$ and
	$h_{x_j}=\th^{-1}(v_{x_j}-\ell_{x_j}v), (j=1,2,\cdots,n)$. From the
	symmetry condition of $a^{jk}$, we see that
	\begin{equation}\label{a1} \sum_{j,k=1}^n
		a^{jk}(\ell_{x_j}v_{x_j}+\ell_{x_k}v_{x_j})=2\sum_{j,k=1}^n
		a^{jk}\ell_{x_j}v_{x_k}.
	\end{equation}
	From \eqref{a1}, we find that
\begin{eqnarray}\label{a2}
	&& \th \sum_{j,k=1}^n( a^{jk}h_{x_j})_{x_k}=\th\sum_{j,k=1}^n\[ \th^{-1}a^{jk}(v_{x_j}-\ell_{x_j}v)\]_{x_k}\nonumber\\
	&& =\sum_{j,k=1}^n\[a^{jk}(v_{x_j}-\ell_{x_j}v)\]_{x_k}-\sum_{j,k=1}^n a^{jk}(v_{x_j}-\ell_{x_j}v)\ell_{x_k}\\
	&& =\sum_{j,k=1}^n \Big\{
	(a^{jk}v_{x_j})_{x_k}-2a^{jk}\ell_{x_j}v_{x_k}+\big[a^{jk}\ell_{x_j}\ell_{x_k}-(a^{jk}\ell_{x_j})_{x_k}\big]v\Big\}.\nonumber
\end{eqnarray}
 Recall that $\ds\cA= -\sum_{j,k=1}^n a^{jk}\ell_{x_j}\ell_{x_k} -\ell_t $ and put
	\begin{equation}\label{a3}
		\left\{
		\ba{ll}
		\ds I_1\= -\sum_{j,k=1}^n (a^{jk}v_{x_j})_{x_k}+\cA v,\\
		\ns\ds I_2\= dv+\sum_{j,k=1}^n \[2a^{jk}\ell_{x_j}v_{x_k}+(a^{jk}\ell_{x_j})_{x_k}v\]dt.
		\ea
		\right.
	\end{equation}
	Then, by $dh=\th^{-1}(dv-\ell_t v)dt$ and (\ref{a2}), we get
	\begin{equation}\label{a4}
		\th\[ dh-\sum_{j,k=1}^n (a^{jk}h_{x_j})_{x_k}dt\] =I_1dt+I_2.
	\end{equation}
	Consequently,
	\begin{equation}\label{a5} 2f\th \[ -\sum_{j,k=1}^n( a^{jk}v_{x_j})_{x_k}+\cA v\] \[
		du-\sum_{j,k=1}^n \big(a^{jk}u_{x_j}\big)_{x_k}dt\] =2f I_1^2dt+2fI_1I_2.
	\end{equation}

	\ss
	
	\par {\bf Step 2. } In this step, we compute  $2fI_1I_2$. Noting that
	\begin{equation}\label{a6} \sum_{j,k,j',k'=1}^n
		\big(a^{jk}a^{j'k'}\ell_{x_j'}v_{x_j}v_{x_k}\big)_{x_{k'}}=-\sum_{j,k,j',k'=1}^n
		\big(a^{j'k'}a^{jk}\ell_{x_j}v_{x_{j'}}v_{x_{k'}}\big)_{x_{k}},
	\end{equation}
	we get
	\begin{equation}\label{a9} \ba{ll}
		\ds 4f\sum_{j,k=1}^n a^{jk}\ell_{x_j}v_{x_k}\[ -\sum_{j,k=1}^n (a^{jk} v_{x_j})_{x_k}+\cA v\]\\
		\ns\ds =-4\sum_{j,k=1}^n (fa^{jk}a^{j'k'} \ell_{x_{j'}}v_{x_j}v_{x_{k'}})_{x_k}+4\sum_{j,k=1}^n  a^{jk}( fa^{j'k'} \ell_{x_{j'}})_{x_k}v_{x_j}v_{x_{k'}}\\
		\ns\ds \q + 4\sum_{j,k=1}^n f a^{jk}a^{j'k'}\ell_{x_{j'}}v_{x_j}v_{x_kx_{k'}}+2\cA f \sum_{j,k=1}^n a^{jk} \ell_{x_j} (v^2)_{x_k}\\
		\ns\ds =-2\sum_{j,k=1}^n \[ \sum_{j',k'=1}^n \(  2 fa^{jk}a^{j'k'}\ell_{x_{j'}}v_{x_j}v_{x_{k'}}-f a^{jk}a^{j'k'}\ell_{x_j} v_{x_{j'}}v_{x_{k'}}\) - f\cA a^{jk}\ell_{x_j}v^2\]_{x_k}\\
		\ns\ds \q +2\sum_{j,k,j',k'=1}^n \[ 2 a^{jk'} \big( fa^{j'k} \ell_{x_{x_j}}\big)_{x_{k'}} -\big( f\a^{jk} a^{j'k'}\ell_{x_{j'}}\big)_{x_{k'}}\] v_{x_j}v_{x_k}-2\sum_{j,k=1}^n ( f\cA a^{jk} \ell_{x_j})_{x_k}v^2,
		\ea
	\end{equation}
	and
	\begin{equation}\label{b1}
		\ba{ll}
		\ds 2f\sum_{j,k=1}^n (a^{jk}\ell_{x_j})_{x_k}v\[ -\sum_{j,k=1}^n (a^{jk}v_{x_j})_{x_k}+\cA v\]\\
		\ns\ds =-2\sum_{j,k,j',k'=1}^n \[ fa^{jk}(a^{j'k'}\ell_{x_{j'}})_{x_{k'}}vv_{x_j}\]_{x_k}+2 f\sum_{j,k,j',k'=1}^n(a^{j'k'}\ell_{x_{j'}})_{x_{k'}}a^{jk}v_{x_j}v_{x_k}\\
		\ns\ds\q  +2\sum_{j,k,j'k'=1}^n \big[a^{jk}f(a^{j'k'}\ell_{x_{j'}})_{x_{k'}}\big]_{x_j}v_{x_k}v +2f\sum_{j,k=1}^n (a^{jk}\ell_{x_j})_{x_k}\cA v^2.
		\ea
	\end{equation}
	Using It\^o's formula, we have
	\begin{eqnarray}\label{a10}
		&& 2f\[ -\sum_{j,k=1}^n (a^{jk}v_{x_j})_{x_k} +\cA v\] dv\nonumber\\
		&& =-2\sum_{j,k=1}^n ( fa^{jk}v_{x_j}dv)_{x_k}+2 f\sum_{j,k=1}^n a^{jk}v_{x_j}d v_{x_k}+2\sum_{j,k=1}^n a^{jk}f_{x_j}v_{x_k}dv+2 f\cA vdv\\
		&& =-2\sum_{j,k=1}^n (fa^{jk}v_{x_j}dv)_{x_k}+d\( f \sum_{j,k=1}^n a^{jk}v_{x_j}v_{x_k}+f\cA v^2\)-\sum_{j,k=1}^n (fa^{jk})_tv_{x_j}v_{x_k}dt \nonumber\\
		&& \q -(f\cA)_t v^2dt-f\sum_{j,k=1}^n a^{jk}d v_{x_j}dv_{x_k}-f\cA
		(dv)^2+2\sum_{j,k=1}^n a^{jk}f_{x_j}v_{x_k}dv. \nonumber
	\end{eqnarray}

	It follows from (\ref{a9}), (\ref{b1})  and (\ref{a10}) that
	\begin{eqnarray}\label{a11}
		&&\3n\3n 2fI_1I_2\nonumber \\
		&&\3n\3n =-2\sum_{j,k=1}^n \Big\{ \sum_{j',k'=1}^n \[  2
		fa^{jk}a^{j'k'}\ell_{x_{j'}}v_{x_j}v_{x_{k'}}-f
		a^{jk}a^{j'k'}\ell_{x_j} v_{x_{j'}}v_{x_{k'}}\nonumber
		\\
		&& \qq\qq\qq\qq +fa^{jk}(a^{j'k'}\ell_{x_{j'}})_{x_{k'}}vv_{x_j}\]- f\cA a^{jk}\ell_{x_j}v^2\Big\}_{x_k}dt\nonumber\\
		&&   -2\sum_{j,k=1}^n (fa^{jk}v_{x_j}dv)_{x_k}dt+d\( f \sum_{j,k=1}^n a^{jk}v_{x_j}v_{x_k}+f\cA v^2\)\\
		&&  +2\!\!\!\sum_{j,k,j',k'=1}^n\!\! \[ 2 a^{jk'} ( fa^{j'k}
		\ell_{{x_j}})_{x_{k'}}\!\! -( f\a^{jk}
		a^{j'k'}\ell_{x_{j'}})_{x_{k'}}\!\!+f(a^{j'k'}\ell_{x_{j'}})_{x_{k'}}a^{jk}\!
		-\frac{1}{2} (fa^{jk})_t\] v_{x_j}v_{x_k}dt\nonumber
		\\
		&&  +2\sum_{j,k,j'k'=1}^n\big[
		a^{jk}f(a^{j'k'}\ell_{x_{j'}})_{x_{k'}}\big]_{x_j}v_{x_k}vdt
		+\[-2\sum_{j,k=1}^na^{jk}(f\cA)_{x_j}\ell_{x_k}-(f\cA)_t\]v^2dt \nonumber\\
		&&  -f\sum_{j,k=1}^n a^{jk}d v_{x_j}dv_{x_k}-f\cA
		(dv)^2+2\sum_{j,k=1}^n a^{jk}f_{x_j}v_{x_k}dv.  \nonumber
	\end{eqnarray}
	\ss
	
	\par {\bf Step 3. } Combining (\ref{a5}), (\ref{a11}), and noting
	\begin{equation*}\label{a13}
		f\sum_{j,k=1}^n a^{jk}dv_{x_j}dv_{x_k} +f\cA (dv)^2=f\th^2 \sum_{j,k=1}^n a^{jk} [(dh_{x_j}+\ell_{x_j}dh)(dh_{x_k}+\ell_{x_k}dh)]+f\th^2\cA (dh)^2,
	\end{equation*}
	we immediately obtain the desired equality.
	\endpf

	\section{Carleman estimate for stochastic parabolic
		equations}\label{Sec3}

	\par Without loss of generality, in what follows, we
	assume that $0\in G$ and $x_0=0$. Let $t_0\in
	(0,T)$. For $r\in (0, \min_{x\in\overline
		G}|x|_{\dbR^n})$ and $\d_0\in (0,t_0)$, we set
	$$
	\cB_r\=\big\{x\in G\big|\, |x|_{\dbR^n}\leq
	r\big\},\qq  Q_{r,\d_0} \= \cB_r\times
	(t_0-\d_0,t_0+\d_0).
	$$

	To prove Theorem \ref{pa-sucp-th1}, we first
	establish a Carleman estimate by virtue of Lemma
	\ref{c1t1}. For  simplicity, we denote
	\begin{equation}\label{f1}
		\ba{ll}
		\ds A(t,x)\= (a^{jk}(t,x))_{1\le j,k\le n},\q A_0(t)\= (a^{jk}(t,0))_{1\le j,k\le n}\=(a_0^{jk})_{1\le j,k\le n}, \\
		\ns\ds x\=(x_1,x_2,\cdots, x_n)^\top,\qq\q
		A_0(t)^{-1}\=(b^{jk}(t))_{1\le j,k\le n}. \ea
	\end{equation}
	For a fixed number $\mu\geq 1$ to be
	chosen later, define
	\begin{equation}\label{pa-sucp-eq1}
		\begin{cases}\ds
			\si(x,t)=\sqrt{A_0^{-1}x\cd x},\\
			\ns\ds
			\f(s)=s\exp\(\int_{0}^{s}\frac{e^{-\mu\tau}-1}{\tau}d\tau\),\\
			\ns\ds   w(x,t)=\f(\si(x,t)),\q
			\phi(s)\= e^{\mu s}=\frac{\f(s)}{s\f'(s)}.
		\end{cases}
	\end{equation}

	\begin{remark}
		The weight function we use here is the one people used to establish
		the SUCP for deterministic parabolic equations. However, the proof
		of the Carleman estimate (Lemma \ref{pa-sucp-lm1} below) is not a
		trivial generalization of the deterministic ones. In the stochastic
		setting,  some extra terms involving the covariation processes of
		solutions would appear. One needs to handle these terms carefully.
	\end{remark}

	\begin{lemma}\label{pa-sucp-lm1}
		There exist  $r_0=r_0((a^{jk})_{1\leq j,k\leq
			n})>0$, $s_1\in (0,1)$ and   $\l_0,\mu_0>0$ such that $\mu=\mu_0$ and
		for any $\e_0\in (0,s_1r_0)$, there is a constant
		$C>0$ independent on $\e_0$  so that for all $\l\geq\l_0$ and
		$$
		\begin{array}{ll}\ds
			z\in \cH_{r_0,\d_0} \=\big\{z\in
			L^2_\dbF(\Om;C_0([t_0-\d_0,t_0+\d_0];L^2(\cB_{r_0})))\cap
			L^2_\dbF(t_0-\d_0,t_0+\d_0;H_0^1(\cB_{r_0}))|\\
			\ns\ds \qq\qq\qq z=0 \mbox{ in
			}(t_0-\d_0,t_0+\d_0)\times [\cB_{\e_0}\cup
			(\cB_{r_0}\setminus \cB_{s_1r_0})]\big\},
		\end{array}
		$$
		which solves
		\begin{equation}\label{1.2-eq4}
			dz - \sum_{j,k=1}^n (a^{jk}z_{x_j})_{x_k} dt = g_1dt + g_2dW(t) \q
			\mbox{ in }\  Q_{r_0,\d_0}
		\end{equation}
		for some $g_1\in
		L^2_\dbF(t_0-\d_0,t_0+\d_0;L^2(\cB_{r_0}))$ and
		$g_2\in
		L^2_\dbF(t_0-\d_0,t_0+\d_0;W^{1,\infty}(\cB_{r_0}))$,
		the following inequality holds:
		\begin{equation}\label{pa-sucp-lm1-eq1}
			\begin{array}{ll}\ds
				\q\mE\int_{Q_{r_0,\d_0}}\big(\l w^{1-2\l}|\nabla
				z|^2 +
				\l^3w^{-1-2\l}|z|^2\big)dxdt \\
				\ns\ds \leq C\mE\int_{Q_{r_0,\d_0}}w^{2-2\l}
				(g_1^2 + \l^2 w^{-2}g_2^2 + |\nabla g_2|^2)dxdt,
			\end{array}
		\end{equation}
		where $C$ depend on $(a^{jk})_{1\le j,k\le n}, r_0,\d_0$.
	\end{lemma}

	{\it Proof}\,: The proof is long. We divide it into four steps.
	
	\ss
	
{\bf  Step 1. }Let $\ell(t,x) = -\l\ln w(t,x)$, $f=\si^2\phi $ and
$h=z$ in \eqref{c1e2a}. Integrating (\ref{c1e2a}) on $Q_{r_0,\d_0}$
and taking mathematical expectation, we have that
\begin{eqnarray}\label{pa-sucp-eq3}
&&\displaystyle\3n 
2 \mE\int_{Q_{r_0,\d_0}}\! f\th\[-\!\sum_{j,k=1}^n\! (a^{jk}v_{x_j})_{x_k}+\cA v\]
\[dz-\!\sum_{j,k=1}^n\!(a^{jk}z_{x_j})_{x_k}dt\]dx +2 \mE\!\int_{Q_{r_0,\d_0}}\!\sum_{j,k=1}^n\!
(fa^{jk}v_{x_j}dv)_{x_k}dx \nonumber \\
&&\ds + 2\mE\int_{Q_{r_0,\d_0}}
\sum_{j,k=1}^n\Big\{\sum_{j',k'=1}^n\[2fa^{jk}
a^{j'k'}\ell_{x_{j'}}v_{x_j}v_{x_{k'}}-fa^{jk}a^{j'k'}\ell_{x_j}v_{x_{j'}}v_{x_{k'}} +fa^{jk}(a^{j'k'}\ell_{x_{j'}})_{x_{k'}}vv_{x_j}\]\nonumber
\\ &&\qq\qq\qq\qq  -f\cA a^{jk}\ell_{x_j} v^2\Big\}_{x_k}dxdt\nonumber\\
&&\ds   - \mE\int_{Q_{r_0,\d_0}}d\[f\cA
v^2+f\sum_{j,k=1}^n
a^{jk}v_{x_j}v_{x_k}\]dx \\
&&\ds\3n =2 \mE\int_{Q_{r_0,\d_0}}\sum_{j,k=1}^n
c^{jk}v_{x_j}v_{x_k}dxdt
+ \mE\int_{Q_{r_0,\d_0}}\cB v^2dxdt \nonumber\\
&&\ds  +2 \mE\int_{Q_{r_0,\d_0}}f\[-\sum_{j,k=1}^n \big(a^{jk}v_{x_j}\big)_{x_k}+\cA v\]^2dxdt\nonumber\\
&&\ds  - \mE\int_{Q_{r_0,\d_0}}f\th^2\sum_{j,k=1}^n a^{jk}(
dz_{x_j} + \ell_{x_j}dz)(dz_{x_k}+\ell_{x_k}dz)dx-
\mE\int_{Q_{r_0,\d_0}}f\th^2 \cA (dz)^2dx\nonumber\\
\ns&&\ds  +2\mE \int_{Q_{r_0,\d_0}} \sum_{j,k=1}^n
a^{jk}f_{x_j}v_{x_k}dvdx +2\mE\int_{Q_{r_0,\d_0}}
\sum_{j,k,j',k'=1}^n
a^{jk}f_{x_j}v_{x_k}(a^{j'k'}\ell_{x_{j'}})_{x_{k'}}vdxdt,\nonumber
\end{eqnarray}
where $c^{jk}$, $\cA$ and $\cB$ are given by
(\ref{c1e15}).

	Clearly, noting that $w\sim O(\si)$ as $\si\to 0$, 
	\begin{eqnarray}\label{1.2-eq3}
		&&
		2 \mE\int_{Q_{r_0,\d_0}}f\th\[-\sum_{j,k=1}^n
		(a^{jk}v_{x_j})_{x_k}+\cA
		v\]\[dz-\sum_{j,k=1}^n(a^{jk}z_{x_j})_{x_k}dt\]dx\nonumber
		\\
		&& = 2
		\mE\int_{Q_{r_0,\d_0}}f\th\[-\sum_{j,k=1}^n
		(a^{jk}v_{x_j})_{x_k}+\cA v\]g_1dxdt \\
		&& \leq
		\mE\int_{Q_{r_0,\d_0}}f\[-\sum_{j,k=1}^n
		(a^{jk}v_{x_j})_{x_k}+\cA v\]^2dxdt +C
		\mE\int_{Q_{r_0,\d_0}}\th^2w^2g_1^2dxdt.\nonumber
	\end{eqnarray}
	Noting that $z\in \cH_{r_0,\d_0}$, we find that
	\begin{eqnarray}\label{1.2-eq5}
		&&\3n\3n\3n\3n
		2 \mE\int_{Q_{r_0,\d_0}}\sum_{j,k=1}^n (fa^{jk}v_{x_j}dv)_{x_k}dx \nonumber\\
		&&\3n\3n\3n\3n + 2
		\mE\!\int_{Q_{r_0,\d_0}}\!\sum_{j,k=1}^n\!\Big\{\sum_{j',k'=1}^n\!\big[2fa^{jk}
		a^{j'k'}\ell_{x_{j'}}v_{x_j}v_{x_{k'}}\!-fa^{jk}a^{j'k'}\ell_{x_j}v_{x_{j'}}v_{x_{k'}}\!+fa^{jk}(a^{j'k'}\ell_{x_{j'}})_{x_{k'}}vv_{x_j}
		\big]\nonumber
		\\
		&&\qq\qq\q  -f\cA a^{jk}\ell_{x_j} v^2\Big\}_{x_k}dxdt=0.
	\end{eqnarray}
	\par Recalling that  $\ell(t,x) = -\l\ln w(t,x)$ and $\si=\sqrt{A_0^{-1}x\cd x}$,  it is
	easy to see that
	\begin{equation}\label{si1}
		\ba{ll}
		\ds \si_{x_j}=\frac{1}{\si}\sum_{k=1}^n b^{jk}x_k,\q \n \si=\frac{A_0^{-1}x^\top}{\si},\\
		\ns\ds \si_{x_jx_k}=\frac{b^{jk}}{\si}
		-\frac{\si_{x_j}\si_{x_k}}{\si},\q  \si_t =\frac{\pa_tA_{0}^{-1}x\cd
			x}{2\sqrt{A_0^{-1}x\cd x}}, \ea
	\end{equation}
	and that
	\begin{equation}\label{1.2-eq2}
		\ba{ll}
		\ds \ell_{t}=-\l w^{-1}w_{t}=-\l\frac{\si_t}{\si\phi},\q  \ell_{x_j}=-\l\frac{\si_{x_j}}{\si\phi },\\
		\ns\ds \ell_{x_jx_k}=-\l
		\frac{b^{jk}}{\si^2\phi}+\l\frac{\si_{x_j}\si_{x_k}}{\si^2\phi}(2+\mu\si).
		\ea
	\end{equation}
	It follows from (\ref{si1}) and (\ref{1.2-eq2})  that
	\begin{equation}\label{1.2-eq2-1}
		\Big|f\sum_{j,k=1}^n a^{jk}\mE(dz_{x_j} +
		\ell_{x_j}dz)(dz_{x_k}+\ell_{x_k}dz)- f\cA \mE(dz)^2\Big|\leq
		C\big(\l^2 \mE g_2^2+ w^2\mE|\nabla g_2|^2\big).
	\end{equation}

	Combining  (\ref{pa-sucp-eq3})--(\ref{1.2-eq5}) and
	\eqref{1.2-eq2-1}, we get that
	\begin{eqnarray}\label{t1}
		&& 2 \mE\int_{Q_{r_0,\d_0}}\sum_{j,k=1}^n c^{jk}v_{x_j}v_{x_k}dxdt +
		\mE\int_{Q_{r_0,\d_0}}\cB v^2dxdt  +2\mE \int_{Q_{r_0,\d_0}}
		\sum_{j,k=1}^n a^{jk}f_{x_j}v_{x_k}dvdx \nonumber
		\\
		&& +2\mE\!\int_{Q_{r_0,\d_0}}\! \sum_{j,k,j',k'=1}^n\!\big[
		a^{jk}f(a^{j'k'}\ell_{x_{j'}})_{x_{k'}}\big]_{x_j}v_{x_k}vdxdt\! + \!2
		\mE\!\int_{Q_{r_0,\d_0}}\!\! f\[-\sum_{j,k=1}^n\!
		\big(a^{jk}v_{x_j}\big)_{x_k}\!\! +\!\cA v\]^2dxdt \nonumber
		\\
		&& \le C\dbE \int_{Q_{r_0,\d_0}}w^{2-2\l}(g_1^2+\l^2w^{-2}g_2^2+|\n
		g_2|^2)dxdt.
	\end{eqnarray}

	{\bf Step 2.} Now, we deal with the left hand side of  (\ref{t1}).
	Recalling (\ref{c1e15}) for $c^{jk} $, we have
	\begin{equation}\label{25ep-sucp7}
		\begin{array}{ll}\ds
			\sum_{j,k=1}^n c^{jk}v_{x_k}v_{x_j} =\sum_{j,k,j'k'=1}^n\[
			2a^{jk'}(fa^{j'k}\ell_{x_j'})_{x_{k'}}-(fa^{jk})_{x_{k'}}a^{j'k'}\ell_{x_{j'}}-\frac{1}{2}(fa^{jk})_t\]
			v_{x_j}v_{x_k}.
		\end{array}
	\end{equation}
	For the first and second term in the righthand side of
	\eqref{25ep-sucp7}, we have respectively that
	\begin{equation}\label{t2}
		\ba{ll}
		\ds 2\sum_{j,k,j',k'=1}^n a^{jk'}(fa^{j'k}\ell_{x_{j}})_{x_{k'}}v_{x_j}v_{x_k}\\
		\ns\ds=-2\l\sum_{j,k,j',k'=1}^n a^{jk'}\(\si^2\phi a^{j'k}\frac{\si_{x_j}}{\si\phi}\)_{x_{k'}}v_{x_j}v_{x_k}\\
		\ns\ds =-2\l \sum_{j,k,j',k'=1}^n a^{jk'}( a^{j'k}\si\si_{x_{j'}})_{x_{k'}}v_{x_j}v_{x_k}\\
		\ns\ds =-2\l \sum_{j,k,j',k'=1}^n\[ a^{jk'}a^{j'k}(\si_{x_{j'}}
		\si_{x_{k'}}+\si\si_{x_{j'}x_{k'}})v_{x_j}v_{x_k}+a^{jk'}a^{j'k}_{x_{k'}}\si\si_{x_{j'}}v_{x_j}v_{x_k}\]
		\\
		\ns\ds = -2\l \sum_{j,k,j',k'=1}^n\( a^{jk'}a^{j'k}b^{j'k'}v_{x_j}v_{x_k} +a^{jk'} a^{j'k}_{x_{k'}}\si\si_{x_{j'}}v_{x_j}v_{x_k}\)
		\\
		\ns\ds =-2\l A A_0^{-1}A \n v\cd \n v-2\l  \sum_{j,k,j',k'=1}^n
		a^{jk'}a^{j'k}_{x_{k'}}\si\si_{x_{j'}}v_{x_j}v_{x_k}  \ea
	\end{equation}
	and
	\begin{eqnarray}\label{t3}
		&& -\sum_{j,k,j',k'=1}^n (fa^{jk})_{x_{k'}}a^{j'k'}\ell_{x_{j'}}v_{x_j}v_{x_k}\nonumber\\
		&& = -\sum_{j,k,j',k'=1}^n (\si^2\phi a^{jk})_{x_{k'}}a^{j'k'}\(-\l \frac{\si_{x_{j'}}}{\si\phi }\)v_{x_j}v_{x_k}\\
		&& = \l\sum_{j,k,j',k'=1}^n \(2 a^{jk}a^{j'k'} \si_{x_{j'}}\si_{x_{k'}}  v_{x_j}v_{x_k}+\l\mu \si a^{jk}a^{j'k'} \si_{x_{j'}}\si_{x_{k'}} v_{x_j}v_{x_{k}} +a_{x_{k'}}^{jk}a^{j'k'}\si\si_{x_j}v_{x_j}v_{x_k}\) \nonumber\\
		&& =\l(2+\mu \si) \frac{A_0^{-1}AA_{0}^{-1}x\cd x}{\si^2}(A \n
		v\cd \n v)+\l\sum_{j,k,j',k'=1}^n a^{jk}_{x_{k'}}a^{j'k'}\si\si_{x_j}v_{x_j}v_{x_k}.\nonumber
	\end{eqnarray}
	Recalling $\si=\sqrt{A_0^{-1}x\cd x}$, we find that
	\begin{equation*}\label{t4}
		\frac{A_0^{-1}AA_{0}^{-1}x\cd x}{\si^2} I-AA_0^{-1}=\frac{A_0^{-1} (A-A_0)A_0^{-1}x\cd x}{A_0^{-1}x\cd x} I-(A-A_0)A_0^{-1}.
	\end{equation*}
	This, together  with $|A-A_0|\le C\si $, implies that
	\begin{equation}\label{t6}
		2\l \[\frac{A_0^{-1}AA_0^{-1}x\cd x}{\si^2}(A\n v\cd\n v )- AA_0^{-1}A\n\cd \n v\]\le C\l \si A\n v\cd \n v.
	\end{equation}
	Combining (\ref{25ep-sucp7})--(\ref{t6}), we find that there exists
	a constant $\mu_1>0$ such that for any $\mu\ge \mu_1$,
	\begin{equation}\label{t7}
		2\dbE \int_{Q_{r_0,\d_0}}\sum_{j,k=1}^n c^{jk}v_{x_j}v_{x_k}dxdt\ge
		C\l\dbE  \int_{Q_{r_0,\d_0}} \si |\n v|^2dxdt.
	\end{equation}

	Further, by the first equality in \eqref{c1e15}, we have that
	\begin{equation*}\label{23lw}
		\begin{array}{ll}\ds
			\cA\3n&\ds
			=-\sum_{j,k=1}^na^{jk}\ell_{x_j}\ell_{x_k}-\ell_t\\
			\ns&\ds= -\l^2 \frac{1}{\si^2\phi^2}\sum_{j,k=1}^n a^{jk}\si_{x_j}\si_{x_k} -\l \frac{\si_t}{\si\phi }.
		\end{array}
	\end{equation*}
	This, together with the second equality in \eqref{c1e15}, implies
	that
	\begin{equation*}\label{t8}
		\ba{ll}
		\ds \cB=-2\sum_{j,k=1}^n a^{jk}(f\cA)_{x_j}\ell_{x_k}-(f\cA)_t\\
		\ns\ds \q = 2\l^3\frac{1}{\si\phi }A \n \(\frac{1}{\phi }A\n \si \cd \n \si\)\cd \n \si-O(\l^2)\\
		\ns\ds\q \ge C\l^3\mu \frac{1}{\si\phi^2} (A\n\si\cd\n\si)^2-O(\l^3)\frac{1}{\si\phi }.
		\ea
	\end{equation*}
	Thus, there exists a constant $\mu_2>0$ such that for any $\mu\ge
	\mu_2$, we have
	\begin{equation}\label{t10} \dbE \int_{Q_{r_0,\d_0}}\cB
		v^2dxdt\ge C\l^3\dbE\int_{Q_{r_0,\d_0}} w^{-1}v^2dxdt.
	\end{equation}
	Moveover, it is clear that
	\begin{equation}\label{t11} \ba{ll}
		\ds 2\mE\int_{Q_{r_0,\d_0}} \sum_{j,k,j',k'=1}^n
		\big[a^{jk}f(a^{j'k'}\ell_{x_{j'}})_{x_{k'}}\big]_{x_j}v_{x_k}vdxdt\\
		\ns\ds =2\l \dbE\int_{Q_{r_0,\d_0}} \sum_{j,k,j',k'=1}^n \[a^{jk} a^{j'k'}_{x_{k'}}\si\si_{x_{j'}} -a^{jjk}b^{j'k'}a^{j'k'}+(2+\mu\si) \si_{x_{j'}}\si_{x_{k'}} a^{jk}a^{j'k'}\]_{x_j}v_{x_k}vdxdt\\
		\ns\ds \ge -C\dbE \int_{Q_{r_0,\d_0}} w|\n v|^2dxdt-C\l^2\dbE \int_{Q_{r_0,\d_0}} w^{-1}v^2dxdt.
		\ea
	\end{equation}

	From (\ref{t1}) and (\ref{t7})--(\ref{t11}), by setting
	$\mu=\mu_0=\max\{\mu_1,\mu_2\}$, we find that there $\l_0>0$ such
	that for any $\l\ge \l_0$, it holds that
	\begin{equation}\label{t12}
		\ba{ll}
		\ds 2\dbE\int_{Q_{r_0,\d_0}}\sum_{j,k,j',k'=1}^n a^{jk}f_{x_j}v_{x_k}dvdx+\l\int_{Q_{r_0,\d_0}}\(\l^2w^{-1}v^2+w|\n v|^2\)dxdt\\
		\ns\ds + \mE\int_{Q_{r_0,\d_0}}f\[-\sum_{j,k=1}^n \big(a^{jk}v_{x_j}\big)_{x_k}+\cA v\]^2dxdt\\
		\ns\ds \le C\dbE \int_{Q_{r_0,\d_0}}w^{2-2\l}
		\big(g_1^2+\l^2w^{-2}g_2^2+|\n g_2|^2\big)dxdt. \ea
	\end{equation}

	{\bf Step 3.} Now we deal with the first term on the left hand side
	of (\ref{t12}). Clearly,
	\begin{equation}\label{TRR1}
		\ba{ll} \ds 2\dbE\int_{Q_{r_0,\d_0}} \sum_{j,k,j',k'=1}^n
		a^{jk}f_{x_j}v_{x_k}dvdx =2\dbE
		\int_{Q_{r_0,\d_0}}\sum_{j,k,j',k'=1}^n
		a^{jk}(\si\phi+\mu\si^2\phi)\si_{x_j}v_{x_k}dvdx. \ea
	\end{equation}
	From (\ref{a3}), (\ref{a4}), \eqref{1.2-eq4}, we get that
	\begin{equation}\label{tt1} \ba{ll}
		\ds \th \[ du-\sum_{j,k=1}^n (a^{jk}u_{x_j})_{x_k}\]dt\\
		\ns\ds =\[ -\sum_{j,k=1}^n (a^{jk}v_{x_j})_{x_k}+\cA v\]dt +dv+\sum_{j,k=1}^n \[ 2a^{jk}\ell_{x_j}v_{x_k} +(a^{jk}\ell_{x_j})_{x_k}v\]dt\\
		\ns\ds =\th \(g_1dt+g_2dW(t)\). \ea
	\end{equation}
	Noting that $v=\th u$, we have that
	\begin{equation*}\label{tt2}
		\ba{ll} \ds -\th \sum_{j,k=1}^n (a^{jk}u_{x_j})_{x_k} =
		-\sum_{j,k=1}^n (a^{jk}v_{x_j})_{x_k}+\cA v+\ell_t v+\sum_{j,k=1}^n
		\[ 2a^{jk}\ell_{x_j}v_{x_k} +(a^{jk}\ell_{x_j})_{x_k}v\]. \ea
	\end{equation*}
	Consequently, recalling (\ref{1.2-eq2}), 
	\begin{equation}\label{tt3}
		\ba{ll}
		\ds 2\l\frac{1}{\si\phi }\sum_{j,k=1}^n a^{jk}\si_{x_j}v_{x_k}\\
		\ns\ds =\th \sum_{j,k=1}^n (a^{jk}u_{x_j})_{x_k}-\sum_{j,k=1}^n (a^{jk}v_{x_j})_{x_k}+\cA v+\ell_t v+\sum_{j,k=1}^n (a^{jk}\ell_{x_j})_{x_k}v.
		\ea
	\end{equation}
	Recalling $v=\th u$, we have
	\begin{equation}\label{trr1}
		dv=\th du+\th\ell_tudt=\th\[g_1dt+g_2dW(t)+\sum_{j,k=1}^n
		(a^{jk}u_{x_j})_{x_k}dt+\ell_t udt\].
	\end{equation}
	Combining with (\ref{TRR1}), (\ref{tt3})  and (\ref{trr1}), we
	obtain
	\begin{eqnarray}\label{trr3}
		&& 2\dbE\int_{Q_{r_0,\d_0}} \sum_{j,k,j',k'=1}^n a^{jk}f_{x_j}v_{x_k}dvdx \nonumber\\
		&& =2\dbE \int_{Q_{r_0,\d_0}}\frac{1}{2\l} \si^2\phi^2(1+\mu\si) 2\l\frac{1}{\si\phi}
		\sum_{j,k=1}^n a^{jk}\si_{x_j}v_{x_k}dvdx\\
		&&  =\frac{1}{\l}\dbE \Big[\int_{Q_{r_0,\d_0}}
		\si^2\phi^2(1+\mu\si)\cJ_1\cJ_2dxdt+\int_{Q_{r_0,\d_0}}
		\si^2\phi^2(1+\mu\si)\sum_{j,k=1}^n
		(a^{jk}\ell_{x_j})_{x_k}vdvdx\Big], \nonumber
	\end{eqnarray}
	where
	\begin{equation}\label{trr4}
		\ba{ll}
		\ds \cJ_1=\th \sum_{j,k=1}^n (a^{jk}u_{x_j})_{x_k}-\sum_{j,k=1}^n (a^{jk}v_{x_j})_{x_k}+\cA v+\ell_t v,\\
		\ns\ds \cJ_2 =\th\[g_1+\sum_{j,k=1}^n (a^{jk}u_{x_j})_{x_k}+\ell_t
		u\]. \ea
	\end{equation}
	From \eqref{trr4}, we find that
	\begin{equation}\label{trr5}
		\ba{ll} \ds\frac{1}{\l}\dbE  \int_{Q_{r_0,\d_0}}\si^2\phi^2(1+\mu\si)\cJ_1\cJ_2dxdt\\
		\ns\ds \ge\dbE \int_{Q_{r_0,\d_0}}\frac{1}{\l}
		\si^2\phi^2(1+\mu\si) \Big\{\frac{1}{2}\th^2\[\sum_{j,k=1}^n
		(a^{jk}u_{x_j})_{x_k}\]^2- \[-\sum_{j,k=1}^n (a^{jk}v_{x_j})_{x_k}+\cA v\]^2
		\\
		\ns\ds \qq\qq\qq\qq\qq\qq -Cw^{-2\l}g_1^2-C\l^2 w^{-3-2\l}u^2\Big\}dxdt\\
		\ns\ds \ge \dbE C\l^{-1} \int_{Q_{r_0,\d_0}}\Big\{w^{2} \[-\sum_{j,k=1}^n
		(a^{jk}v_{x_j})_{x_k}+\cA v\]^2-C\(w^{2-2\l}g_1^2+\l^2
		w^{-1-2\l}u^2\)\Big\}dxdt \ea
	\end{equation}
	Recalling $z\in \cH_{r_0,\d_0}$ and using It\^o formula, we get that
	\begin{equation}\label{trr6}
		\ba{ll}
		\ds \frac{1}{\l}\dbE  \int_{Q_{r_0,\d_0}}\si^2\phi^2(1+\mu\si)\sum_{j,k=1}^n (a^{jk}\ell_{x_j})_{x_k}vdvdx\\
		\ns\ds = -\frac{1}{\l }\dbE \int_{Q_{r_0,\d_0}}\si^2\phi^2(1+\mu\si) \[\sum_{j,k=1}^n (a^{jk}\ell_{tx_j})_{x_k}v^2-\sum_{j,k=1}^n (a^{jk}\ell_{x_j})_{x_k}\th^2 g_2^2\]dxdt\\
		\ns\ds \ge -C \int_{Q_{r_0,\d_0}}  w^{-2\l}(g_2^2+ u^2)dxdt.
		\ea
	\end{equation}

	{\bf Step 4. } Combining (\ref{t12}), (\ref{trr3}), (\ref{trr5}) and
	(\ref{trr6}),  recalling that $v=w^{-\l}u$, we get
	\eqref{pa-sucp-lm1-eq1}. \endpf

	
	\section{Proof of Theorem \ref{pa-sucp-th1}}\label{Sec4}
	
	
	Assume that $z \in L^2_\dbF(\Om;C([t_0-\d_0,t_0+\d_0];
	L^2(\cB_{r_0})))\cap L^2_\dbF(t_0-\d_0,t_0+\d_0; H^1(\cB_{r_0}))$
	satisfies that
	\begin{equation}\label{pa-sucp-lm2-eq1}
		dz -\sum_{j,k=1}^n(a^{jk}z_{x_j})_{x_k}dt = \hat
		a\cdot\nabla zdt + \hat b zdt + \hat c z dW(t)
		\q\mbox{ in }\; Q_{r_0,\d_0},
	\end{equation}
	where $\hat a\in
	L^\infty_\dbF(t_0-\d_0,t_0+\d_0;L^\infty(\cB_{r_0};\dbR^n))$,
	$\hat b\in
	L^\infty_\dbF(t_0-\d_0,t_0+\d_0;L^\infty(\cB_{r_0}))$
	and $\hat c\in
	L^\infty_\dbF(t_0-\d_0,t_0+\d_0;W^{1,\infty}(\cB_{r_0}))$. Set
	\begin{equation}\label{yy1}
		M=|\hat
		a|_{L^\infty_\dbF(t_0-\d_0,t_0+\d_0;L^\infty(\cB_{r_0};\dbR^n))}+|\hat
		b|_{L^\infty_\dbF(t_0-\d_0,t_0+\d_0;L^\infty(\cB_{r_0}))}+|\hat
		c|_{L^\infty_\dbF(t_0-\d_0,t_0+\d_0;W^{1,\infty}(\cB_{r_0}))}.
	\end{equation}
	For any fixed constant $0<\d_1<\d_0$, in the rest of this section,
	we denote by $\cC\=\cC(\d_1,\d_0, M,$ $(a^{jk})_{1\le j,k\le n},
	\hat a, \hat b, \hat c,r_0) $  a generic constant which may change
	from line to line.
	\begin{lemma}\label{pa-sucp-lm2}
		
		There exists a constant $\cC
		> 1$ such that for all $0<r_2<r_1<2r_1<r_0$, it holds that
		\begin{equation}\label{pa-sucp-lm2-eq2}
			\ba{ll} \ds |z|_{L^2_\dbF(t_0-\d_1,t_0+\d_1;L^2(\cB_{r_1}))} \3n&\ds
			\leq \cC |z|_{L^2_\dbF(t_0-\d_0,t_0+\d_0;L^2(\cB_{r_2}))}^{\e_0}
			|z|_{L^2_\dbF(t_0-\d_0,t_0+\d_0;L^2(\cB_{r_0}))}^{1-\e_0}\\
			\ns&\ds \q  +e^{\cC\ln \f (2r_0/3)-\ln
				\f(r_2/3)}|z|_{L^2_\dbF(t_0-\d_0,t_0+\d_0;L^2(\cB_{r_0}))}, \ea
		\end{equation}
		where $\f$ given in (\ref{pa-sucp-eq1}) and 
		\begin{equation}\label{pa-sucp-lm2-eq3}
			\e_0=\frac{\ln \f(2r_0/3) -  \ln \f(r_1) }{\ln \f(2r_0/3) - \ln
				\f (r_2/2)}.
		\end{equation}
	\end{lemma}
	\begin{remark}
		Lemma \ref{pa-sucp-lm2} is a generalization of the classical three
		cylinder inequality of parabolic equations (e.g., \cite{Vessella}).
		Besides being an important tool to prove the strong unique
		continuation, as the deterministic case, it can be used to solve
		some inverse problems for stochastic parabolic equations. However,
		this is beyond the scope of this paper and will be presented in our
		future works.
	\end{remark}
	{\it Proof}\,:  Let $t_1\in(0,(\d_0-\d_1)/4)$. Set $T_1= \d_0 -t_1$
	and $T_2=\d_0 - 2t_1$. Let $\psi\in C_0^2(t_0-\d_0,t_0+\d_0)$ such
	that
	\begin{equation}\label{pa-sucp14}
		\psi(t)=\left\{
		\begin{array}{ll}\ds
			0, &\mbox{ if } t\in [t_0-\d_0,t_0-T_1]\cup
			[t_0+T_1,t_0+\d_0],\\
			\ns\ds 1, &\mbox{ if } t\in [t_0-T_2,t_0+T_2],\\
			\ns\ds
			\exp\(-\frac{\d_0^3(T_2-t+t_0)^4}{(T_1-t+t_0)^3t_1^4}\),
			&
			\mbox{ if }t\in (t_0+T_2,t_0+T_1),\\
			\ns\ds
			\exp\(-\frac{\d_0^3(T_2+t-t_0)^4}{(T_1+t-t_0)^3t_1^4}\),
			& \mbox{ if }t\in (t_0-T_1,t_0-T_2).
		\end{array}
		\right.
	\end{equation}

	Let $\a_0 =r_2/2$ and let $f\in C_0^2(0,r_0)$
	such that
	\begin{equation}\label{pa-sucp14.1}
		f(t)=\left\{
		\begin{array}{ll}\ds
			0, &\mbox{ if } t\in [0, \a_0]\cup [3r_0/4, r_0],\\
			\ns\ds 1, &\mbox{ if } t\in
			[3\a_0/2,2r_0/3],
		\end{array}
		\right.
	\end{equation}
	and that
	\begin{equation}\label{pa-sucp14.2}
		\left\{
		\begin{array}{ll}\ds |f'|\leq C_f/\a_0,\q
			|f''|\leq C_f/\a_0^2 &\mbox{
				in }[\a_0,3\a_0/2],\\
			\ns\ds |f'|\leq C_f/r_0,\q |f''|\leq C_f/r_0^2
			&\mbox{ in }[r_0/2,3r_0/4],
		\end{array}
		\right.
	\end{equation}
	where $C_f$ is an absolute constant.
	
	Let us choose $\zeta$ as
	\begin{equation}\label{pa-sucp15}
		\zeta(x,t)=f(|x|)\psi(t),\qq \mbox{ if }(x,t)\in
		Q_{r_0,\d_0}.
	\end{equation}
	Then $\tilde z\=\zeta z$ solves
	\begin{equation}\label{1.2-eq6}
		d\tilde z - \sum_{j,k=1}^n(a^{jk}\tilde z_{x_j})_{x_k}dt =
		(\hat a\cdot\nabla \tilde z  + \hat b \tilde z
		-\tilde f)dt + \hat c \tilde z dW(t) \q\mbox{ in
		}\; Q_{r_0,\d_0},
	\end{equation}
	where $$\tilde f = \zeta_t z + 2\sum_{j,k=1}^n
	a^{jk} z_{x_j}\zeta_{x_k}+ \sum_{j,k=1}^n
	(a^{jk}\zeta_{x_k})_{x_j} z.$$

	By applying the inequality
	\eqref{pa-sucp-lm1-eq1} to $\tilde z$, we obtain
	that
	\begin{equation}\label{pa-sucp13}
		\begin{array}{ll}\ds
			\q\mE\int_{Q_{r_0,\d_0}}\big(\l
			w^{-2\l}|\nabla\tilde z|^2 +
			\l^3w^{-2-2\l}\tilde z^2 \big)dxdt \\
			\ns\ds \leq
			\cC\mE\int_{Q_{r_0,\d_0}}w^{2-2\l}\big(\big|\hat
			a\cdot\nabla \tilde z  + \hat b \tilde z -\tilde
			f\big|^2+|\nabla (\hat c\tilde
			z)|^2\big)dxdt
			+\cC\l^2\mE\int_{Q_{r_0,\d_0}}w^{-2\l}\big|\hat c
			\tilde z\big|^2dxdt.
		\end{array}
	\end{equation}
	Denote by
	$$
	\begin{array}{ll}\ds
		K_1'=\Big\{ (x,t)\in
		\dbR^{n+1}\Big|\,\frac{3}{2}\a_0\leq
		|x|\leq\frac{2r_0}{3},\; t\in
		[t_0-T_1,t_0-T_2]
		\Big\},\\
		\ns\ds K_1''=\Big\{ (x,t)\in
		\dbR^{n+1}\Big|\,\frac{3}{2}\a_0\leq
		|x|\leq\frac{2r_0}{3},\; t\in
		[t_0+T_2,t_0+T_1]
		\Big\},\\
		\ns\ds K_2=\Big\{ (x,t)\in
		\dbR^{n+1}\Big|\,\a_0\leq
		|x|\leq\frac{3\a_0}{2},\; t\in [t_0-T_1,t_0+T_1]
		\Big\},\\
		\ns\ds K_3=\Big\{ (x,t)\in
		\dbR^{n+1}\Big|\,\frac{2r_0}{3} \leq
		|x|\leq\frac{3r_0}{4},\; t\in
		[t_0-T_1,t_0+T_1]
		\Big\},\\
		\ns\ds K_4=\Big\{ (x,t)\in
		\dbR^{n+1}\Big|\,\frac{3}{2}\a_0\leq
		|x|\leq\frac{2r_0}{3},\; t\in
		[t_0-T_2,t_0+T_2]
		\Big\},\\
		\ns\ds K_1= K_1'\cup K_1'',\qq K_5 =
		Q_{r_0,\d_0}\setminus \cup_{i=1}^4 K_i.
	\end{array}
	$$
	Clearly, we have that $Q_{r_0,\d_0}=\cup_{i=1}^5 K_i$ and $K_i\cap
	K_j=\emptyset$ for $i\neq j$, $i,j = 1, 2, 3, 4, 5$.
	
	It follows from \eqref{pa-sucp13} that for every
	$\l\geq \l_0$,
	\begin{equation}\label{pa-sucp16}
		\begin{array}{ll}\ds
			\mE\int_{K_4}\big(\l w^{-2\l}|\nabla\tilde z|^2
			+
			\l^3w^{-2-2\l}\tilde z^2\big)dxdt\\
			\ns\ds \leq
			J_1+J_2+\cC \mE\int_{K_4}w^{2-2\l}(z^2 + |\nabla
			z|^2)dxdt,
		\end{array}
	\end{equation}
	where
	\begin{equation}\label{pa-sucp16-1}
		\begin{array}{ll}\ds
			J_1\3n&\ds\=-\mE\int_{K_1}\big(\l
			w^{-2\l}|\nabla\tilde z|^2 +
			\l^3w^{-2-2\l}\tilde z^2\big)dxdt\\
			\ns&\ds \q + \cC\mE\!\int_{K_1}\!w^{2-2\l}\big(\big|\hat
			a\cdot\nabla \tilde z  + \hat b \tilde z -\tilde f\big|^2\!+|\nabla
			(\hat c\tilde
			z)|^2\big)dxdt+\cC\l^2\mE\!\int_{K_1}\!w^{-2\l}\big|\hat c \tilde
			z\big|^2dxdt,
		\end{array}
	\end{equation}
	and
	\begin{equation}\label{pa-sucp16-2}
		\begin{array}{ll}\ds
			J_2\3n&\ds\=  \cC\mE\int_{K_2\cup
				K_3}w^{2-2\l}\big(\big|\hat a\cdot\nabla \tilde
			z + \hat b \tilde z -\tilde f\big|^2+|\nabla
			(\hat c\tilde
			z)|^2\big)dxdt\\
			\ns&\ds\q+\cC\l^2\mE\int_{K_2\cup
				K_3}w^{-2\l}\big|\hat c \tilde z\big|^2dxdt.
		\end{array}
	\end{equation}
	By  \eqref{pa-sucp16}, we obtain that there is a $\l_1\geq 0$ such
	that for all $\l\geq\max\{\l_0,\l_1\}$,
	\begin{equation}\label{pa-sucp17}
		\mE\int_{K_4}\big[\l w^{-2\l}|\nabla(z
		\zeta)|^2 + \l^3w^{-2-2\l}z^2
		\zeta^2\big]dxdt\leq J_1+J_2.
	\end{equation}

	It follows from \eqref{pa-sucp16-2}  that
	\begin{equation}\label{pa-sucp27}
		J_2\leq  \cC \f(\a_0 )^{2-2\l} \mE\int_{K_2}(|\nabla z|^2
		+ z^2)dxdt + \cC \f\(\frac{3r_0}{4}\)^{2-2\l} \mE\int_{K_3}( |\nabla
		z|^2 +  z^2)dxdt.
	\end{equation}
	Now we estimate $J_1$. Let
	\begin{equation}\label{pa-sucp18}
		\cE(t,x;\l)=\psi(t)^2\[\cC
		w^{2-2\l}+\cC\(\frac{\psi'(t)}{\psi(t)}\)^2 - \l^3w(x,t)^{-2-2\l}\].
	\end{equation}
	From \eqref{pa-sucp16-1}, we get that
	$$
	\begin{array}{ll}\ds
		J_1\leq \mE\int_{K_1}\cE(t,x;\l) z^2 dxdt +\mE\int_{K_1}\big(\cC
		w^{2-2\l}-\l w^{-2\l} \big)|\nabla z|^2.
	\end{array}
	$$
	It follows from the choice of $w$  that there is
	a $\l_2>0$ such that for all $\l\geq\l_2$,
	$$
	\cC w^{2-2\l}-\l w^{-2\l}\leq 0
	$$
	and
	$$
	\l^3w^{-2-2\l}\geq 2\cC w^{2-2\l}.
	$$
	Thus, for all $\l\geq\l_2$,
	\begin{equation}\label{pa-sucp19}
		J_1\leq\mE\int_{K_1}\l^3\psi^2w^{-1-2\l}\[w^3\(\frac{\psi'(t)}{\psi(t)}\)^2\frac{\cC}{\l^3}-1\]z^2dxdt.
	\end{equation}
	We first handle the case that $(t,x)\in K_1'$. The estimate on
	$K_1''$ would be similar. Recalling $T_1-T_2=t_1$, we get that
	\begin{equation}\label{y4}
		\psi'(t)=-\psi(t)\(\frac{\d_0^3}{t_1^4}\)\frac{4(T_2+t-t_0)^3(T_1+t-t_0)-3(T_2+t-t_0)^4}{(T_1+t-t_0)^4}.
	\end{equation}
	Set
	\begin{equation}\label{y1}
		\e(t,x;\l)\= w^3\(\frac{\psi'(t)}{\psi(t) }\)^2\frac{\cC}{\l^3}-1.
	\end{equation}
	Then we know that there exists a constant $\cC_1>0$  such that
	\begin{equation}\label{y2}
		\e(t,x;\l) \le -\frac{1}{2}+\frac{\cC_1w^3}{\l^3(T_1+t-t_0)^8}
		\q\mbox{ on }K_1'.
	\end{equation}
	Set
	$$
	K_{1,\l}\= \Big\{(t,x)\in K_1'\|~
	-\frac{1}{2}+\frac{\cC_1w^3}{\l^3(T_1+t-t_0)^8}\ge 0\Big\}.
	$$
	Obviously, we have
	\begin{equation}\label{y3}
		\ba{ll} \ds\dbE
		\int_{K_1'}\l^3\psi^2w^{-1-2\l}\[w^3\(\frac{\psi'(t)}{\psi(t)}\)^2\frac{\cC}{\l^3}-1\]z^2dxdt\le
		\cC \dbE \int_{K_{1,\l}} \l^3 \psi'^2w^{2-2\l}z^2dxdt. \ea
	\end{equation}
	From (\ref{y4}), we see that $\phi'^{2}/\psi^{3/2}$ is bounded in
	$K_{1,\l}$. This, together with \eqref{y3}, implies that
	\begin{equation}\label{o1}
		\ba{ll} \ds\dbE
		\int_{K_1'}\l^3\psi^2w^{-1-2\l}\[w^3\(\frac{\psi'(t)}{\psi(t)}\)^2
		\frac{\cC}{\l^3}-1\]z^2dxdt\le \cC \dbE \int_{K_{1,\l}}
		\l^3w^{2-2\l}\phi^{3/2}z^2dxdt. \ea
	\end{equation}
	On the other hand, when $(t,x)\in K_{1,\l}$, we have
	\begin{equation}\label{y7}
		\ba{ll}
		\ds \frac{T_1+t-t_0}{\d_0}\le \(\frac{2\cC_1w^3}{\l^3\d_0^8}\)^{1/8}.
		\ea
	\end{equation}
	After choosing $\l \ge \l_1\=
	\[\frac{2^{9}\cC_1w^3}{t_1^8}\]^{1/3}$, we get that
	\begin{equation}\label{y88}
		|T_2+t-t_0|=|T_2-T_1+T_1+t-t_0|\ge t_1-\(\frac{2\cC_1w^3}{\l^3}\)^{1/8}\ge \frac{t_1}{2}.
	\end{equation}
	From \eqref{y88}, we know there exists a constant $\l_3$,
	independent of $r_2$, such that for any $\l\ge \l_3$,
	\begin{equation}\label{y8}
		\!\!\begin{array}{ll}\ds
			\f\(\frac{3r_0}{4}\)^{2\l-2}\!\psi^{3/2}w^{2-2\l}\3n&\ds\!\le\!
			\exp\!\Big\{\!-\!\[\(\frac{\l^3\d_0^8}{2\cC_1w^3}\)^{3/8}\frac{t_1^4}{2^4}\frac{1}{t_1^4}\]^{3/2}
			\!\!+(2\!-\!2\l)\[\ln
			w\!-\!\ln \f\(\frac{3r_0}{4}\)\]\!\Big\}\\
			\ns&\ds\le 1.
		\end{array}
	\end{equation}
	Combining (\ref{o1}) with (\ref{y8}), and doing a similar argument
	on region $K_1''$, we conclude that
	\begin{equation}\label{y9}
		J_1\le \cC \f\(\frac{3r_0}{4}\)^{2-2\l}\int_{K_1} z^2dxdt.
	\end{equation}

	\par Notice that $r_1\in(3\a_0/2,2r_0/3)$ and denote by $K_4^{(r_1)}$ the
	region $\{(x,t)\in K_4||x|\leq r_1\}$. By \eqref{pa-sucp17},
	\eqref{pa-sucp27} and  \eqref{y9}, we obtain that for all
	$\l\geq\max\{\l_3,\l_5\}$,
	\begin{equation}\label{pa-sucp28}
		\begin{array}{ll}\ds
			\q\mE\int_{K_4^{(r_1)}} z^2dxdt\\
			\ns\ds\leq \f(r_1)^{2\l+2}\mE\int_{K_4} z^2
			w^{-2-2\l}dxdt\\
			\ns\ds \leq \cC \f(r_1)^{2\l+2}\[\f
			(\a_0)^{2-2\l}\mE\int_{K_2}(|\nabla z|^2 + z^2)dxdt +
			\f\(\frac{3r_0}{4}\)^{2-2\l}\mE\int_{K_1} z^2dxdt\\
			\ns\ds \hspace{2.6cm} +
			\f\(\frac{3r_0}{4}\)^{2-2\l}\mE\int_{K_3}(|\nabla z|^2 + z^2)dxdt\].
		\end{array}
	\end{equation}

	Let $\zeta_1\in
	C^\infty_0(Q_{2\a_0,\d_0}\setminus
	Q_{\a_0/2,\d_0})$ such that $\zeta_1=1$ in
	$K_2$. By It\^o's formula, we get that
	\begin{equation}\label{1.2-eq7}
		d(\zeta_1^2 z^2)=2\zeta_{1}\zeta_{1,t}z^2 +
		2\zeta_1^2 zdz + \zeta_1^2 (dz)^2.
	\end{equation}
	Integrating \eqref{1.2-eq7} on
	$Q_{2\a_0,\d_0}\setminus Q_{\a_0/2,\d_0}$ and
	taking mathematical expectation, we find that
	\begin{equation}\label{1.2-eq8}
		\begin{array}{ll}\ds
			0 \3n&\ds=2\mE\int_{Q_{2\a_0,\d_0}\setminus
				Q_{\a_0/2,\d_0}}\zeta_{1}\zeta_{1,t}z^2dxdt +
			2\mE\int_{Q_{2\a_0,\d_0}\setminus
				Q_{\a_0/2,\d_0}}\zeta_1^2 zdzdx\\
			\ns&\ds\q + \mE\int_{Q_{2\a_0,\d_0}\setminus
				Q_{\a_0/2,\d_0}}\zeta_1^2 (dz)^2dx.
		\end{array}
	\end{equation}
	It follows from \eqref{pa-sucp-lm2-eq1} that
	\begin{eqnarray}\label{1.2-eq9}
		&& \3n\3n \mE\int_{Q_{2\a_0,\d_0}\setminus
			Q_{\a_0/2,\d_0}}\zeta_1^2
		zdzdx\nonumber\\
		&&\3n\3n\3n\3n = \mE\int_{Q_{2\a_0,\d_0}\setminus
			Q_{\a_0/2,\d_0}}\zeta_1^2
		z\[\sum_{j,k=1}^n(a^{jk}z_{x_j})_{x_k}  + \hat
		a\cdot\nabla z  + \hat b z\]dxdt\\
		&& \3n\3n\3n\3n= \mE\!\int_{Q_{2r_0,\d_0}\setminus
			Q_{\a_0/2,\d_0}}\!\!\!\zeta_1^2\! \sum_{j,k=1}^n\!
		a^{jk}z_{x_j}z_{x_k}dxdt \!+\! 2\mE\!\int_{Q_{2\a_0,\d_0}\setminus
			Q_{\a_0/2,\d_0}}\!\!\zeta_1 \!\sum_{j,k=1}^n\!
		a^{jk}z_{x_j}\zeta_{1,x_k} z
		dxdt\nonumber\\
		&& \3n\3n + \mE\int_{Q_{2\a_0,\d_0}\setminus
			Q_{\a_0/2,\d_0}} \zeta_1^2 z \big(a\cdot\nabla z
		+ \hat b z\big)dxdt.\nonumber
	\end{eqnarray}
	Combing \eqref{1.2-eq8} and \eqref{1.2-eq9}, we
	find that for any $\e>0$,
	$$
	\begin{array}{ll}\ds
		\q\mE\int_{Q_{2\a_0,\d_0}\setminus
			Q_{\a_0/2,\d_0}}\zeta_1^2
		\sum_{j,k=1}^n a^{jk}z_{x_j}z_{x_k}dxdt\\
		\ns\ds \leq \e\mE\int_{Q_{2\a_0,\d_0}\setminus
			Q_{\a_0/2,\d_0}}\zeta_1^2 |\nabla z|^2 dxdt +
		\frac{\cC}{\e}\mE\int_{Q_{2\a_0,\d_0}\setminus
			Q_{\a_0/2,\d_0}}z^2dxdt.
	\end{array}
	$$
	This, together with \eqref{1.2-eq1}, implies
	that
	\begin{equation}\label{1.2-eq10}
		\begin{array}{ll}\ds
			\q s_0\mE\int_{Q_{2\a_0,\d_0}\setminus
				Q_{\a_0/2,\d_0}}\zeta_1^2
			|\nabla z|^2dxdt\\
			\ns\ds \leq \e\mE\int_{Q_{2\a_0,\d_0}\setminus
				Q_{\a_0/2,\d_0}}\zeta_1^2 |\nabla z|^2 dxdt +
			\frac{\cC}{\e}\mE\int_{Q_{2\a_0,\d_0}\setminus
				Q_{\a_0/2,\d_0}}z^2dxdt.
		\end{array}
	\end{equation}
	By choosing $\e = s_0/2$, from \eqref{1.2-eq10}
	and the definition of $\zeta_1$, we obtain that
	\begin{equation}\label{pa-sucp29}
		\mE\int_{K_2} |\nabla z|^2 dxdt \leq \cC
		\mE\int_{Q_{2\a_0,\d_0}\setminus
			Q_{\a_0/2,\d_0}}z^2dxdt.
	\end{equation}
	Similarly, we can get that
	\begin{equation}\label{pa-sucp30}
		\mE\int_{K_3} |\nabla z|^2 dxdt \leq \cC
		\mE\int_{Q_{r_0,\d_0}\setminus Q_{r_0/3,\d_0}}
		z^2 dxdt.
	\end{equation}

	Put
	$$
	\eta=|z|_{L^2_\dbF(t_0-\d_0,t_0+\d_0;L^2(\cB_{r_2}))}, \q
	\tilde\eta=|z|_{L^2_\dbF(t_0-\d_0,t_0+\d_0;L^2(\cB_{r_0}))}.
	$$

	Adding $\mE\int_{Q_{r_1,\d_1}} z^2 dxdt$ to both sides of
	\eqref{pa-sucp28}, from \eqref{pa-sucp29} and \eqref{pa-sucp30}, we
	obtain that
	\begin{equation}\label{pa-sucp31}
		|z|_{L^2_\dbF(t_0-\d_1,t_0+\d_1;L^2(\cB_{r_1}))}^2
		\leq \cC\[ \(\frac{\f(3\a_0/2)}{\f(r_1)}\)^{2-2\l }\eta^2 +
		\(\frac{\f(2r_0/3)}{\f(r_1)}\)^{2-2\l} \tilde\eta^2 \].
	\end{equation}

	Set
	\begin{equation}\label{pa-sucp32}
		\l_4 = \frac{\ln \tilde \eta-\ln \eta }{\ln \f(2r_0/3)-\ln
			\f(r_2/2)}.
	\end{equation}
	If $\l_4\ge \max_{i=1,2,3}\{\l_i\}$, then, by choosing
	in \eqref{pa-sucp31} $\l=\l_4$,  we get that
	\begin{equation}\label{pa-sucp33}
		|z|_{L^2_\dbF(t_0-\d_1,t_0+\d_1;L^2(\cB_{r_1}))}
		\leq \cC \eta^{\e_0}\tilde\eta^{1-\e_0},
	\end{equation}
	where $\e_0$ given in (\ref{pa-sucp-lm2-eq3}).
	\par  If $\l_4 <  \max_{i=1,2,3}\{\l_i\}$, then, by
	\eqref{pa-sucp32}, we have that
	\begin{equation}\label{pa-sucp35}
		|z|_{L^2_\dbF(t_0-\d_1,t_0+\d_1;L^2(\cB_{r_1}))}
		\leq \tilde\eta \leq
		e^{\cC\ln \f (2r_0/3)-\ln \f(r_2/3)} \eta.
	\end{equation}
	This, together with  \eqref{pa-sucp33},  yields
	(\ref{pa-sucp-lm2-eq2}).
	\endpf
	
	\vspace{0.1cm}
	
	Now we are in a position to prove Theorem
	\ref{pa-sucp-th1}.
	
	\vspace{0.1cm}
	
	{\it Proof of Theorem \ref{pa-sucp-th1}}\,:
	Recalling that  for any $N\in\dbN$ and $r>0$, it
	holds that
	\begin{equation}\label{pa-sucp-eq12}
		\mE\int_{Q_{r,\d_0}}|y(t,x)|^2dxdt=O(r^{2N}).
	\end{equation}
	Applying Lemma \ref{pa-sucp-lm2} to the equation
	\eqref{pa-sucp0},  by \eqref{pa-sucp-eq12} and
	passing to the limit as $r_2$ tends to $0$,  we
	obtain that
	\begin{equation}\label{pa-sucp-eq13}
		\mE|y|^2_{L^2(Q_{r_1,\d_1})}\leq \cC e^{-\cC N}, \q
		\mbox{for every }N\in\dbN,
	\end{equation}
	where $\cC$ is independent of $N$. Passing to the
	limit as $N\to+\infty$, \eqref{pa-sucp-eq13}
	yields that $y = 0$ in $Q_{r_1,\d_1}$,
	$\dbP$-a.s.  By iteration, it follows that $y =
	0$ in $Q$, $\dbP$-a.s.
	\endpf


	\section{Further comments}
	
	
	As far as we know, Theorem \ref{pa-sucp-th1} is
	the first result concerning the SUCP for
	stochastic PDEs. Compared with the fruitful
	study of the SUCP for deterministic PDEs, lots
	of things should be done and some of them seem
	to be very interesting and challenging.
	
	\begin{itemize}
		\item \textbf{The SUCP for stochastic parabolic
			equations with nonsmooth coefficients.} In
		\cite{Escauriaza,Koch2}, the authors show that
		the SUCP for deterministic parabolic equations
		holds when the coefficients $a$ and $b$ are
		integrable in some weighted spaces. We believe
		these results can be generalized in the
		stochastic setting. However, to this end, one
		has to develop $L^p$ Carleman estimate for
		stochastic parabolic equations. It seems to us
		that this is a fascinating but difficult
		problem.
		
		\item \textbf{SUCP for other type of stochastic
			PDEs.} SUCP is also studied for some other type of
		PDEs, such as wave equations (e.g.
		\cite{Lebeau,Vessella2}) and plate equations
		(e.g. \cite{Tarama}). It is diverting to see
		whether these results hold for corresponding
		stochastic PDEs. However, although some
		Carleman estimates have been obtained for some
		other stochastic PDEs (e.g. \cite{Luqi7,Lu,LZ2,Zhangxu3}), as far as we
		know, they cannot be used to establish SUCP
		for these equations.

		\item \textbf{Applications with the SUCP.} As
		we said in the introduction, there are lots of applications of SUCP
		for deterministic PDEs. It is quite interesting to investigate
		applications of SUCP for stochastic PDEs. Some of them can be done
		easily. For example, our result implies approximate controllability
		of backward stochastic parabolic equations. Another example is that
		following the idea in \cite{Vessella1}, one can get some results for
		some inverse problems of stochastic parabolic equations. Details of
		these two applications are beyond of the scope of this paper and
		will be investigated in our forthcoming papers.
	\end{itemize}

	\section*{Acknowledgement}
	
	We appreciate Professor Luis  Escauriaza for
	pointing out some references for SUCP for
	deterministic parabolic equations.

	

\end{document}